\documentclass{llncs}

\usepackage[utf8]{inputenc}
\usepackage{amssymb}
\usepackage{amsmath}

\usepackage{graphics}
\usepackage{alltt,listings}

\usepackage{tikz,tikz-qtree,tikz-qtree-compat,colortbl}
\usepackage{tikz,tikz-qtree,tikz-qtree-compat,adjustbox}
\usepackage{etoolbox}

\newcommand{\circledgenerator}[2][]{\tikz[baseline=(char.base)]{\node[shape = circle, draw, inner sep = 2pt,fill=black!40](char) {\phantom{\ifblank{#1}{#2}{#1}}};\node at (char.center) {\makebox(0,0){#2}};}}

\newcommand{\circlednongap}[2][]{\tikz[baseline=(char.base)]{\node[shape = circle, draw, inner sep = 2pt,fill=black!20](char) {\phantom{\ifblank{#1}{#2}{#1}}};\node at (char.center) {\makebox(0,0){#2}};}}
\newcommand{\circledgap}[2][]{\tikz[baseline=(char.base)]{\node[shape = circle, draw, inner sep = 2pt](char) {\phantom{\ifblank{#1}{#2}{#1}}};\node at (char.center) {\makebox(0,0){#2}};}}

\newcommand{\circledseed}[2][]{\tikz[baseline=(char.base)]{\node[shape = circle, draw, inner sep = 2pt,fill=black!05](char) {\phantom{\ifblank{#1}{#2}{#1}}};\node at (char.center) {\makebox(0,0){#2}};}}
\newcommand{\circlednonseed}[2][]{\tikz[baseline=(char.base)]{\node[shape = circle, draw, inner sep = 2pt,fill=black!30](char) {\phantom{\ifblank{#1}{#2}{#1}}};\node at (char.center) {\makebox(0,0){#2}};}}

\robustify{\circlednongap}
\robustify{\circledgenerator}
\robustify{\circledgap}
\robustify{\circledseed}
\robustify{\circlednonseed}

\newcommand{\gap}[1]{\circledgap[00]{\large \phantom{#1}}}
\newcommand{\generator}[1]{\circledgenerator[00]{\large #1}}
\newcommand{\nongap}[1]{\circlednongap[00]{\large #1}}
\newcommand{\nonseed}[1]{\circledseed[00]{\large #1}}
\newcommand{\seed}[1]{\circlednonseed[00]{\large #1}}

\usetikzlibrary{turtle}\usetikzlibrary{shapes.geometric}

\title{New Eliahou semigroups and verification of the Wilf conjecture for genus up to 65}
\date{\today}
\author{Maria Bras-Amorós, César Marín Rodríguez}

\institute{Universitat Rovira i Virgili}

\begin{document}
\maketitle

\begin{abstract}
We give a graphical reinterpretation of the seeds algorithm to explore the tree of numerical semigroups. We then exploit the seeds algorithm to find all the Eliahou semigroups of genus up to $65$. Since all these semigroups satisfy the Wilf conjecture, this shows that the Wilf conjecture holds up to genus $65$.
  \end{abstract}

\section{Introduction}
A {\em numerical semigroup} is a cofinite submonoid of ${\mathbb N}_0$.
See \cite{RG} for a general reference on numerical semigroups.
The elements in the complement of a numerical semigroup in ${\mathbb N}_0$ are denoted the {\em gaps} of the semigroup. The {\em genus} of the semigroup is the number of its gaps.

If a numerical semigroup $\Lambda$ is $\{\lambda_0=0<\lambda_1<\dots\}$, define its {\em multiplicity} as $m(\Lambda)=\lambda_1$.
%, its {\em jump} as $u(\Lambda)=\lambda_2-\lambda_1$ and its {\em second jump} as $v(\Lambda)=\lambda_3-\lambda_2$.
Define its
{\em Frobenius number} $F(\Lambda)$ as its largest gap and its 
{\em conductor} $c(\Lambda)$ as its largest gap plus one.
If $c(\Lambda)=\lambda_L$, then the elements $\lambda_0,\lambda_1,\dots,\lambda_{L-1}$ are called the {\em left elements} of $\Lambda$.

%A numerical semigroup will be called {\em ordinary} if $c(\Lambda)=m(\Lambda)=\lambda_1$,
%{\em pseudo-ordinary} if $c(\Lambda)=m(\Lambda)+u(\Lambda)=\lambda_2$ and {\em pseudo-pseudo-ordinary} if $c(\Lambda)=m(\Lambda)+u(\Lambda)+v(\Lambda)=\lambda_3$.

An element $\lambda_s\geq c(\Lambda)$ is an order-$i$ {\em seed} of $\Lambda$ if
$\lambda_s+\lambda_i\neq \lambda_j+\lambda_k$ for all $i<j\leq k<s$.
The {\em right primitive elements} of a numerical semigroup are its order-$0$ seeds.
In general, the primitive elements (or minimal generators) of a numerical semigroup are those elements of the semigroup that can not be obtained as a sum of two smaller semigroup elements.

In 1978, Hebert S. Wilf conjectured that
for any numerical semigroup with conductor $c$, with $L$ left elements and with set of primitive elements equal to $P$,
it holds $c\leq L\cdot\#P$ \cite{wilf}. More than forty years later the conjecture is still open. It has been verified for all semigroups of genus up to $60$ by Jean Fromentin and Florent Hivert \cite{FromentinHivert}.
An important step to approach the Wilf conjecture is a sufficient condition found by Shalom Eliahou \cite{eliahou}. Semigroups not satisfying it are very unusual. We denote them {\em Eliahou semigroups}. 

If we take away a primitive element from a numerical semigroup we obtain another semigroup with genus increased by one.
We can organize all numerical semigroups in an infinite tree rooted at ${\mathbb N}_0$ and such that the children of a node are the semigroups obtained taking away one by one its right primitive elements.
In Figure~\ref{f:tree} one can see the lowest genus semigroups organized in the semigroup tree. Each semigroup is represented by its non-gaps which are either colored with dark gray if they are right primitive elements or with light gray if they are not.

\begin{figure}
  \resizebox{1.000000\textwidth}{!}{\begin{tikzpicture}[sibling distance=1.75em,grow'=right, every node/.style = {align=left}]\tikzset{level 1/.style={level distance=3cm}}\tikzset{level 2/.style={level distance=5cm}}\tikzset{level 3/.style={level distance=6cm}}\tikzset{level 4/.style={level distance=8cm}}\tikzset{level 5+/.style={level distance=10cm}}\tikzset{level 7+/.style={level distance=12cm}}\Tree[.{$\nongap{0}\generator{1} \dots$} [.{$\nongap{0}\gap{1}\generator{2} \generator{3} \dots$} [.{$\nongap{0}\gap{1}\gap{2}\generator{3} \generator{4} \generator{5} \dots$} [.{$\nongap{0}\gap{1}\gap{2}\gap{3}\generator{4} \generator{5} \generator{6} \generator{7} \dots$} [.{$\nongap{0}\gap{1}\gap{2}\gap{3}\gap{4}\generator{5} \generator{6} \generator{7} \generator{8} \generator{9} \dots$} [.{$\nongap{0}\gap{1}\gap{2}\gap{3}\gap{4}\gap{5}\generator{6} \generator{7} \generator{8} \generator{9} \generator{10} \generator{11} \dots$}  ]
[.{$\nongap{0}\gap{1}\gap{2}\gap{3}\gap{4}\nongap{5}\gap{6}\generator{7} \generator{8} \generator{9} \nongap{10}\generator{11} \dots$}  ]
[.{$\nongap{0}\gap{1}\gap{2}\gap{3}\gap{4}\nongap{5}\nongap{6}\gap{7}\generator{8} \generator{9} \nongap{10}\nongap{11}\dots$}  ]
[.{$\nongap{0}\gap{1}\gap{2}\gap{3}\gap{4}\nongap{5}\nongap{6}\nongap{7}\gap{8}\generator{9} \nongap{10}\nongap{11}\dots$}  ]
[.{$\nongap{0}\gap{1}\gap{2}\gap{3}\gap{4}\nongap{5}\nongap{6}\nongap{7}\nongap{8}\gap{9}\nongap{10} \nongap{11}\dots$}  ]
 ]
[.{$\nongap{0}\gap{1}\gap{2}\gap{3}\nongap{4}\gap{5}\generator{6} \generator{7} \nongap{8}\generator{9} \dots$} [.{$\nongap{0}\gap{1}\gap{2}\gap{3}\nongap{4}\gap{5}\gap{6}\generator{7} \nongap{8}\generator{9} \generator{10} \nongap{11}\dots$}  ]
[.{$\nongap{0}\gap{1}\gap{2}\gap{3}\nongap{4}\gap{5}\nongap{6}\gap{7}\nongap{8} \generator{9} \nongap{10}\generator{11} \dots$}  ]
[.{$\nongap{0}\gap{1}\gap{2}\gap{3}\nongap{4}\gap{5}\nongap{6}\nongap{7}\nongap{8}\gap{9}\nongap{10} \nongap{11}\dots$}  ]
 ]
[.{$\nongap{0}\gap{1}\gap{2}\gap{3}\nongap{4}\nongap{5}\gap{6}\generator{7} \nongap{8}\nongap{9}\dots$} [.{$\nongap{0}\gap{1}\gap{2}\gap{3}\nongap{4}\nongap{5}\gap{6}\gap{7}\nongap{8} \nongap{9}\nongap{10}\generator{11} \dots$}  ]
 ]
[.{$\nongap{0}\gap{1}\gap{2}\gap{3}\nongap{4}\nongap{5}\nongap{6}\gap{7}\nongap{8} \nongap{9}\dots$}  ]
 ]
[.{$\nongap{0}\gap{1}\gap{2}\nongap{3}\gap{4}\generator{5} \nongap{6}\generator{7} \dots$} [.{$\nongap{0}\gap{1}\gap{2}\nongap{3}\gap{4}\gap{5}\nongap{6} \generator{7} \generator{8} \nongap{9}\dots$} [.{$\nongap{0}\gap{1}\gap{2}\nongap{3}\gap{4}\gap{5}\nongap{6}\gap{7}\generator{8} \nongap{9}\generator{10} \nongap{11}\dots$}  ]
[.{$\nongap{0}\gap{1}\gap{2}\nongap{3}\gap{4}\gap{5}\nongap{6}\nongap{7}\gap{8}\nongap{9} \nongap{10}\generator{11} \dots$}  ]
 ]
[.{$\nongap{0}\gap{1}\gap{2}\nongap{3}\gap{4}\nongap{5}\nongap{6}\gap{7}\nongap{8} \nongap{9}\dots$}  ]
 ]
[.{$\nongap{0}\gap{1}\gap{2}\nongap{3}\nongap{4}\gap{5}\nongap{6} \nongap{7}\dots$}  ]
 ]
[.{$\nongap{0}\gap{1}\nongap{2}\gap{3}\nongap{4} \generator{5} \dots$} [.{$\nongap{0}\gap{1}\nongap{2}\gap{3}\nongap{4}\gap{5}\nongap{6} \generator{7} \dots$} [.{$\nongap{0}\gap{1}\nongap{2}\gap{3}\nongap{4}\gap{5}\nongap{6}\gap{7}\nongap{8} \generator{9} \dots$} [.{$\nongap{0}\gap{1}\nongap{2}\gap{3}\nongap{4}\gap{5}\nongap{6}\gap{7}\nongap{8}\gap{9}\nongap{10} \generator{11} \dots$}  ]
 ]
 ]
 ]
 ]
 ] \end{tikzpicture}}

\caption{Lowest depth nodes of the semigroup tree.}
\label{f:tree}
\end{figure}

In Section~\ref{s:revisit} we recall the {\em seeds algorithm} (\cite{seeds}) to explore the semigroup tree. In Section~\ref{s:example} we give a graphical explanation of the algorithm by running it over a particular example.
In Section~\ref{s:EliahouWilf} we give the list of Eliahou semigroups output by our algorithm. This list allows to verify the Wilf conjecture for genus up to $65$.

\section{The bitstream of gaps and the bitstream of seeds of a numerical semigroup}
\label{s:revisit}

A bitstream is a finite sequence $a=a_0\dots a_\ell$ where $a_i$ is either $0$ or $1$ for every $i$.
For our purposes, we can indistinctly use $a$ for $a_0\dots a_\ell$ and for any bitstream of the form $a_0\dots a_\ell\underbrace{0\dots0}_{k}$ for any positive integer $k$.

Suppose a semigroup $\Lambda$ has conductor $c$. We encode its gaps as the bitstream $$G(\Lambda)=g_0\dots g_{c-1}$$ with $g_i=0$ if $i+1\in \Lambda$ is a gap and $g_i=1$ otherwise.
We encode its seeds as the bitstream $$S(\Lambda)=s_0\dots s_{c-1}$$ with $s_i=1$ if $i-\lambda_j$ is an order-$j$ seed of $\Lambda$ where $j$ is the unique non-negative integer such that $\lambda_j\leq i <\lambda_{j+1}$.

%\begin{example}
%For instance, the semigroup $\Lambda=\{0,3,6,8,9,10,\dots\}$, with conductor $c=8$ has $G(\Lambda)=11011010$ and $S(\Lambda)=10110111$.
%\end{example}

In \cite{seeds} we presented an algorithm to explore the tree of numerical semigroups by recursively computing the bitstream of gaps and the bitstream of seeds of a numerical semigroup from those of its parent.

In Figure~\ref{f:treeseeds} one can see the sequences $G$ and $S$ for the lowest genus semigroups organized in the semigroup tree. To make it easier to read, we represented each $1$ in the sequences with a dark circle with its position written inside, and each $0$ in the sequence with a light gray circle with its position also written inside.

\begin{figure}
  \resizebox{1.000000\textwidth}{!}{\begin{tikzpicture}[sibling distance=3em,grow'=right, every node/.style = {align=left}]\tikzset{level 1/.style={level distance=3cm}}\tikzset{level 2/.style={level distance=5cm}}\tikzset{level 3/.style={level distance=6cm}}\tikzset{level 4/.style={level distance=8cm}}\tikzset{level 5+/.style={level distance=10cm}}\tikzset{level 7+/.style={level distance=12cm}}\Tree[.{$\begin{array}{l}\nonseed{0}\\\seed{0} \end{array}$} [.{$\begin{array}{l}\seed{0}\nonseed{1}\\\seed{0} \seed{1} \end{array}$} [.{$\begin{array}{l}\seed{0}\seed{1}\nonseed{2}\\\seed{0} \seed{1} \seed{2} \end{array}$} [.{$\begin{array}{l}\seed{0}\seed{1}\seed{2}\nonseed{3}\\\seed{0} \seed{1} \seed{2} \seed{3} \end{array}$} [.{$\begin{array}{l}\seed{0}\seed{1}\seed{2}\seed{3}\nonseed{4}\\\seed{0} \seed{1} \seed{2} \seed{3} \seed{4} \end{array}$} [.{$\begin{array}{l}\seed{0}\seed{1}\seed{2}\seed{3}\seed{4}\nonseed{5}\\\seed{0} \seed{1} \seed{2} \seed{3} \seed{4} \seed{5} \end{array}$}  ]
[.{$\begin{array}{l}\seed{0}\seed{1}\seed{2}\seed{3}\nonseed{4}\seed{5}\nonseed{6}\\\seed{0} \seed{1} \seed{2} \nonseed{3} \seed{4} \seed{5} \seed{6} \end{array}$}  ]
[.{$\begin{array}{l}\seed{0}\seed{1}\seed{2}\seed{3}\nonseed{4}\nonseed{5}\seed{6}\nonseed{7}\\\seed{0} \seed{1} \nonseed{2} \nonseed{3} \nonseed{4} \seed{5} \seed{6} \seed{7} \end{array}$}  ]
[.{$\begin{array}{l}\seed{0}\seed{1}\seed{2}\seed{3}\nonseed{4}\nonseed{5}\nonseed{6}\seed{7}\nonseed{8}\\\seed{0} \nonseed{1} \nonseed{2} \nonseed{3} \nonseed{4} \nonseed{5} \seed{6} \seed{7} \seed{8} \end{array}$}  ]
[.{$\begin{array}{l}\seed{0}\seed{1}\seed{2}\seed{3}\nonseed{4}\nonseed{5}\nonseed{6}\nonseed{7}\seed{8}\nonseed{9}\\\nonseed{0} \nonseed{1} \nonseed{2} \nonseed{3} \nonseed{4} \nonseed{5} \nonseed{6} \seed{7} \seed{8} \seed{9} \end{array}$}  ]
 ]
[.{$\begin{array}{l}\seed{0}\seed{1}\seed{2}\nonseed{3}\seed{4}\nonseed{5}\\\seed{0} \seed{1} \nonseed{2} \seed{3} \seed{4} \seed{5} \end{array}$} [.{$\begin{array}{l}\seed{0}\seed{1}\seed{2}\nonseed{3}\seed{4}\seed{5}\nonseed{6}\\\seed{0} \nonseed{1} \seed{2} \seed{3} \seed{4} \seed{5} \seed{6} \end{array}$}  ]
[.{$\begin{array}{l}\seed{0}\seed{1}\seed{2}\nonseed{3}\seed{4}\nonseed{5}\seed{6}\nonseed{7}\\\nonseed{0} \seed{1} \nonseed{2} \seed{3} \nonseed{4} \seed{5} \seed{6} \seed{7} \end{array}$}  ]
[.{$\begin{array}{l}\seed{0}\seed{1}\seed{2}\nonseed{3}\seed{4}\nonseed{5}\nonseed{6}\nonseed{7}\seed{8}\nonseed{9}\\\nonseed{0} \nonseed{1} \nonseed{2} \nonseed{3} \nonseed{4} \nonseed{5} \nonseed{6} \seed{7} \seed{8} \seed{9} \end{array}$}  ]
 ]
[.{$\begin{array}{l}\seed{0}\seed{1}\seed{2}\nonseed{3}\nonseed{4}\seed{5}\nonseed{6}\\\seed{0} \nonseed{1} \nonseed{2} \nonseed{3} \seed{4} \seed{5} \seed{6} \end{array}$} [.{$\begin{array}{l}\seed{0}\seed{1}\seed{2}\nonseed{3}\nonseed{4}\seed{5}\seed{6}\nonseed{7}\\\nonseed{0} \nonseed{1} \nonseed{2} \seed{3} \seed{4} \seed{5} \seed{6} \seed{7} \end{array}$}  ]
 ]
[.{$\begin{array}{l}\seed{0}\seed{1}\seed{2}\nonseed{3}\nonseed{4}\nonseed{5}\seed{6}\nonseed{7}\\\nonseed{0} \nonseed{1} \nonseed{2} \nonseed{3} \nonseed{4} \seed{5} \seed{6} \seed{7} \end{array}$}  ]
 ]
[.{$\begin{array}{l}\seed{0}\seed{1}\nonseed{2}\seed{3}\nonseed{4}\\\seed{0} \nonseed{1} \seed{2} \seed{3} \seed{4} \end{array}$} [.{$\begin{array}{l}\seed{0}\seed{1}\nonseed{2}\seed{3}\seed{4}\nonseed{5}\\\nonseed{0} \seed{1} \seed{2} \seed{3} \seed{4} \seed{5} \end{array}$} [.{$\begin{array}{l}\seed{0}\seed{1}\nonseed{2}\seed{3}\seed{4}\nonseed{5}\seed{6}\nonseed{7}\\\seed{0} \nonseed{1} \seed{2} \seed{3} \nonseed{4} \seed{5} \seed{6} \seed{7} \end{array}$}  ]
[.{$\begin{array}{l}\seed{0}\seed{1}\nonseed{2}\seed{3}\seed{4}\nonseed{5}\nonseed{6}\seed{7}\nonseed{8}\\\nonseed{0} \nonseed{1} \seed{2} \nonseed{3} \nonseed{4} \nonseed{5} \seed{6} \seed{7} \seed{8} \end{array}$}  ]
 ]
[.{$\begin{array}{l}\seed{0}\seed{1}\nonseed{2}\seed{3}\nonseed{4}\nonseed{5}\seed{6}\nonseed{7}\\\nonseed{0} \nonseed{1} \nonseed{2} \nonseed{3} \nonseed{4} \seed{5} \seed{6} \seed{7} \end{array}$}  ]
 ]
[.{$\begin{array}{l}\seed{0}\seed{1}\nonseed{2}\nonseed{3}\seed{4}\nonseed{5}\\\nonseed{0} \nonseed{1} \nonseed{2} \seed{3} \seed{4} \seed{5} \end{array}$}  ]
 ]
[.{$\begin{array}{l}\seed{0}\nonseed{1}\seed{2}\nonseed{3}\\\nonseed{0} \seed{1} \seed{2} \seed{3} \end{array}$} [.{$\begin{array}{l}\seed{0}\nonseed{1}\seed{2}\nonseed{3}\seed{4}\nonseed{5}\\\nonseed{0} \seed{1} \nonseed{2} \seed{3} \seed{4} \seed{5} \end{array}$} [.{$\begin{array}{l}\seed{0}\nonseed{1}\seed{2}\nonseed{3}\seed{4}\nonseed{5}\seed{6}\nonseed{7}\\\nonseed{0} \seed{1} \nonseed{2} \seed{3} \nonseed{4} \seed{5} \seed{6} \seed{7} \end{array}$} [.{$\begin{array}{l}\seed{0}\nonseed{1}\seed{2}\nonseed{3}\seed{4}\nonseed{5}\seed{6}\nonseed{7}\seed{8}\nonseed{9}\\\nonseed{0} \seed{1} \nonseed{2} \seed{3} \nonseed{4} \seed{5} \nonseed{6} \seed{7} \seed{8} \seed{9} \end{array}$}  ]
 ]
 ]
 ]
 ]
 ] \end{tikzpicture}}
\caption{Bitstream of gaps and bitstream of seeds of the semigroups in the lowest depth nodes of the semigroup tree.}
\label{f:treeseeds}
\end{figure}

The updating algorithm is based on the next results, which are proved in \cite{seeds}.
Suppose that $\lambda_s$ is a right primitive element of $\Lambda$ (hence, $s\geq L$) and let $\tilde\Lambda=\Lambda\setminus\{\lambda_s\}$.

\begin{enumerate}
\item {\bf Old-order recycled seeds:}

Suppose $i<L$. Any order-$i$ seed $\lambda_t$ of $\Lambda$ with \,$t>s$ is also 
an order-$i$ seed of $\tilde\Lambda$.

\item {\bf Old-order new seeds:}
Suppose $i<L$.
Then,
$\lambda_t>\lambda_s$,
with $\lambda_t$ not an order-$i$ seed of $\Lambda$,
is an order-$i$ seed of $\tilde\Lambda$ 
if and only if either
\begin{itemize}
%\item[{\rm (1)}] $\lambda_t$ is an order-$i$ seed of $\Lambda$ (hence, an old-order recycled seed)
\item $i<L-1$, \,$\lambda_t=\lambda_s+\lambda_{i+1}-\lambda_i$\, and
\,$\lambda_s$ is an order-$(i+1)$ seed of $\Lambda$
\item $i=L-1$, $\lambda_s=c$, and either $\left\{\begin{array}{l}\lambda_t=\lambda_s+\lambda_L-\lambda_{L-1}\\\lambda_t=\lambda_s+\lambda_L-\lambda_{L-1}+1\end{array}\right.$
\item $i=L-1$, $\lambda_s=c+1$, and \,$\lambda_t=\lambda_s+\lambda_L-\lambda_{L-1}$
\end{itemize}

\item{\bf New-order seeds:}
Suppose $i\geq L$. 
Then,
\begin{itemize}
\item If \,$i<s-2$, then $\tilde\Lambda$ has no order-$i$ seeds.
\item If \,$i=s-2$, then the only order-$i$ seed of \,$\tilde\Lambda$ is \,$\lambda_s+1$.
\item If \,$i=s-1$, then the only order-$i$ seeds of \,$\tilde\Lambda$ are \,$\lambda_s+1$\, and \,$\lambda_s+2$.
\end{itemize}
  \end{enumerate}

\section{A graphical explanation of the algorithm by an example}
\label{s:example}

For the algorithm in \cite{seeds}
the bitstream of seeds is splitted in a table 
with $L$ rows, indexed from $0$ to $L-1$, with the $i$th row containing
$s_{(\lambda_i)},s_{({\lambda_i+1})}\dots,s_{(\lambda_{i+1}-1)}$.

Now we are going to graphically explain the algorithm with an example.
Consider the numerical semigroup $\Lambda=\{0, 8, 16, 18, 19, 24, 26, 27, 30,\dots\}$. Its table of seeds is as follows, where seeds are repesented by black boxes and non-seeds are represented by white boxes.
Notice that its conductor is $30$ and it has three right primitive elements which are $30, 31, 33$.

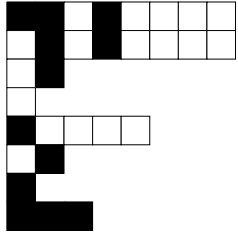
\usetikzlibrary{arrows}\newcommand\ambcolor{\cellcolor{black}\color{white}}
\renewcommand{\arraystretch}{.85}
\newcommand\neteja[2]{\begin{scope}[xshift=#1cm,yshift=-#2cm]\draw[fill=white,draw=none,turtle={home,right,forward,right,forward,right,forward,right,forward}];\end{scope}\begin{scope}[xshift=#1cm,yshift=-#2cm]\draw[fill=gray!25,draw=none,turtle={home,right,forward,right,forward,right,forward,right,forward}];\end{scope}}
\newcommand\seeddos[2]{\begin{scope}[xshift=#1cm,yshift=-#2cm]\filldraw[fill=black,turtle={home,right,forward,right,forward,right,forward,right,forward}];\end{scope}}
\newcommand\nonseeddos[2]{\begin{scope}[xshift=#1cm,yshift=-#2cm]\filldraw[fill=white,turtle={home,right,forward,right,forward,right,forward,right,forward}];\end{scope}}
\newcommand\transparentseed[2]{\begin{scope}[xshift=#1cm,yshift=-#2cm]\filldraw[fill=black,fill opacity=0.4,turtle={home,right,forward,right,forward,right,forward,right,forward}];\end{scope}}
\newcommand\transparentnonseed[2]{\begin{scope}[xshift=#1cm,yshift=-#2cm]\filldraw[fill=white,fill opacity=0.4,turtle={home,right,forward,right,forward,right,forward,right,forward}];\end{scope}}
\newcommand\netejaseed[2]{\neteja{#1}{#2}\seeddos{#1}{#2}}
\newcommand\netejanonseed[2]{\neteja{#1}{#2}\nonseeddos{#1}{#2}}
\newcommand\netejatransparentseed[2]{\neteja{#1}{#2}\transparentseed{#1}{#2}}
\newcommand\netejatransparentnonseed[2]{\neteja{#1}{#2}\transparentnonseed{#1}{#2}}
\begin{center}\scalebox{0.380000}{\begin{tikzpicture}[turtle/distance=1cm]\filldraw [thick,white,fill=gray!25,turtle={home,right,forward,left,right,right,left,forward=8,right,forward,forward,right,forward=6,left,forward,right,forward=1,left,forward,left,forward=4,right,forward,right,forward=3,left,forward,right,forward=1,left,forward,left,forward=2,right,forward,right,forward=3,right,forward=8}];\seeddos{1.000000}{0}\seeddos{2.000000}{0}\nonseeddos{3.000000}{0}\seeddos{4.000000}{0}\nonseeddos{5.000000}{0}\nonseeddos{6.000000}{0}\nonseeddos{7.000000}{0}\nonseeddos{8.000000}{0}\nonseeddos{1.000000}{1}\seeddos{2.000000}{1}\nonseeddos{3.000000}{1}\seeddos{4.000000}{1}\nonseeddos{5.000000}{1}\nonseeddos{6.000000}{1}\nonseeddos{7.000000}{1}\nonseeddos{8.000000}{1}\nonseeddos{1.000000}{2}\seeddos{2.000000}{2}\nonseeddos{1.000000}{3}\seeddos{1.000000}{4}\nonseeddos{2.000000}{4}\nonseeddos{3.000000}{4}\nonseeddos{4.000000}{4}\nonseeddos{5.000000}{4}\nonseeddos{1.000000}{5}\seeddos{2.000000}{5}\seeddos{1.000000}{6}\seeddos{1.000000}{7}\seeddos{2.000000}{7}\seeddos{3.000000}{7}\end{tikzpicture}
}\end{center}Suppose we want to take away the generator c+0=30+0=30.\begin{center}\scalebox{0.380000}
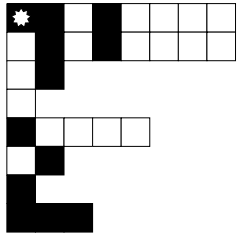
{\begin{tikzpicture}[turtle/distance=1cm]\filldraw [thick,white,fill=gray!25,turtle={home,right,forward,left,right,right,left,forward=8,right,forward,forward,right,forward=6,left,forward,right,forward=1,left,forward,left,forward=4,right,forward,right,forward=3,left,forward,right,forward=1,left,forward,left,forward=2,right,forward,right,forward=3,right,forward=8}];\seeddos{1.000000}{0}\seeddos{2.000000}{0}\nonseeddos{3.000000}{0}\seeddos{4.000000}{0}\nonseeddos{5.000000}{0}\nonseeddos{6.000000}{0}\nonseeddos{7.000000}{0}\nonseeddos{8.000000}{0}\nonseeddos{1.000000}{1}\seeddos{2.000000}{1}\nonseeddos{3.000000}{1}\seeddos{4.000000}{1}\nonseeddos{5.000000}{1}\nonseeddos{6.000000}{1}\nonseeddos{7.000000}{1}\nonseeddos{8.000000}{1}\nonseeddos{1.000000}{2}\seeddos{2.000000}{2}\nonseeddos{1.000000}{3}\seeddos{1.000000}{4}\nonseeddos{2.000000}{4}\nonseeddos{3.000000}{4}\nonseeddos{4.000000}{4}\nonseeddos{5.000000}{4}\nonseeddos{1.000000}{5}\seeddos{2.000000}{5}\seeddos{1.000000}{6}\seeddos{1.000000}{7}\seeddos{2.000000}{7}\seeddos{3.000000}{7}\node[draw, fill=white, star, star points=9, minimum size=7mm] at (1.500000,-.5){};
\end{tikzpicture}
}\end{center}Draw the contour of the new table of seeds.

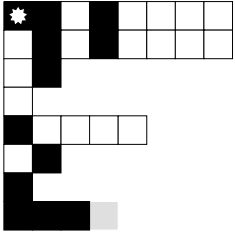
\begin{center}\scalebox{0.380000}{\begin{tikzpicture}[turtle/distance=1cm]\filldraw [thick,white,fill=gray!25,turtle={home,right,forward,left,right,right,left,forward=8,right,forward,forward,right,forward=6,left,forward,right,forward=1,left,forward,left,forward=4,right,forward,right,forward=3,left,forward,right,forward=1,left,forward,left,forward=3,right,forward,right,forward=4,right,forward=8}];\seeddos{1.000000}{0}\seeddos{2.000000}{0}\nonseeddos{3.000000}{0}\seeddos{4.000000}{0}\nonseeddos{5.000000}{0}\nonseeddos{6.000000}{0}\nonseeddos{7.000000}{0}\nonseeddos{8.000000}{0}\nonseeddos{1.000000}{1}\seeddos{2.000000}{1}\nonseeddos{3.000000}{1}\seeddos{4.000000}{1}\nonseeddos{5.000000}{1}\nonseeddos{6.000000}{1}\nonseeddos{7.000000}{1}\nonseeddos{8.000000}{1}\nonseeddos{1.000000}{2}\seeddos{2.000000}{2}\nonseeddos{1.000000}{3}\seeddos{1.000000}{4}\nonseeddos{2.000000}{4}\nonseeddos{3.000000}{4}\nonseeddos{4.000000}{4}\nonseeddos{5.000000}{4}\nonseeddos{1.000000}{5}\seeddos{2.000000}{5}\seeddos{1.000000}{6}\seeddos{1.000000}{7}\seeddos{2.000000}{7}\seeddos{3.000000}{7}\node[draw, fill=white, star, star points=9, minimum size=7mm] at (1.500000,-.5){};
\end{tikzpicture}
}\end{center}Move the old values to the left of the table and fix the old-order recycled seeds.

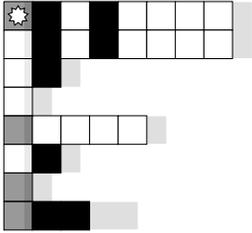
\begin{center}\scalebox{0.380000}{\begin{tikzpicture}[turtle/distance=1cm]\filldraw [thick,white,fill=gray!25,turtle={home,right,forward,left,right,right,left,forward=8,right,forward,forward,right,forward=6,left,forward,right,forward=1,left,forward,left,forward=4,right,forward,right,forward=3,left,forward,right,forward=1,left,forward,left,forward=3,right,forward,right,forward=4,right,forward=8}];\transparentseed{0.300000}{0}\seeddos{1.300000}{0}\nonseeddos{2.300000}{0}\seeddos{3.300000}{0}\nonseeddos{4.300000}{0}\nonseeddos{5.300000}{0}\nonseeddos{6.300000}{0}\nonseeddos{7.300000}{0}\transparentnonseed{0.300000}{1}\seeddos{1.300000}{1}\nonseeddos{2.300000}{1}\seeddos{3.300000}{1}\nonseeddos{4.300000}{1}\nonseeddos{5.300000}{1}\nonseeddos{6.300000}{1}\nonseeddos{7.300000}{1}\transparentnonseed{0.300000}{2}\seeddos{1.300000}{2}\transparentnonseed{0.300000}{3}\transparentseed{0.300000}{4}\nonseeddos{1.300000}{4}\nonseeddos{2.300000}{4}\nonseeddos{3.300000}{4}\nonseeddos{4.300000}{4}\transparentnonseed{0.300000}{5}\seeddos{1.300000}{5}\transparentseed{0.300000}{6}\transparentseed{0.300000}{7}\seeddos{1.300000}{7}\seeddos{2.300000}{7}\node[draw, fill=white, star, star points=9, minimum size=7mm] at (0.800000,-.5){};
\end{tikzpicture}
}\end{center}Obtain the old-order new seeds.

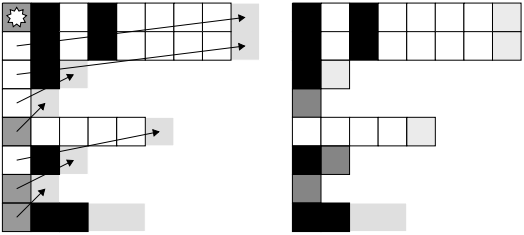
\begin{center}\scalebox{0.380000}{\begin{tikzpicture}[turtle/distance=1cm]\filldraw [thick,white,fill=gray!25,turtle={home,right,forward,left,right,right,left,forward=8,right,forward,forward,right,forward=6,left,forward,right,forward=1,left,forward,left,forward=4,right,forward,right,forward=3,left,forward,right,forward=1,left,forward,left,forward=3,right,forward,right,forward=4,right,forward=8}];\transparentseed{0.000000}{0}\seeddos{1.000000}{0}\nonseeddos{2.000000}{0}\seeddos{3.000000}{0}\nonseeddos{4.000000}{0}\nonseeddos{5.000000}{0}\nonseeddos{6.000000}{0}\nonseeddos{7.000000}{0}\transparentnonseed{0.000000}{1}\seeddos{1.000000}{1}\nonseeddos{2.000000}{1}\seeddos{3.000000}{1}\nonseeddos{4.000000}{1}\nonseeddos{5.000000}{1}\nonseeddos{6.000000}{1}\nonseeddos{7.000000}{1}\transparentnonseed{0.000000}{2}\seeddos{1.000000}{2}\transparentnonseed{0.000000}{3}\transparentseed{0.000000}{4}\nonseeddos{1.000000}{4}\nonseeddos{2.000000}{4}\nonseeddos{3.000000}{4}\nonseeddos{4.000000}{4}\transparentnonseed{0.000000}{5}\seeddos{1.000000}{5}\transparentseed{0.000000}{6}\transparentseed{0.000000}{7}\seeddos{1.000000}{7}\seeddos{2.000000}{7}\node[draw, fill=white, star, star points=9, minimum size=7mm] at (0.500000,-.5){};
\draw[-triangle 60,thick](0.500000,-1.500000) -- (8.500000,-0.500000);
\draw[-triangle 60,thick](0.500000,-2.500000) -- (8.500000,-1.500000);
\draw[-triangle 60,thick](0.500000,-3.500000) -- (2.500000,-2.500000);
\draw[-triangle 60,thick](0.500000,-4.500000) -- (1.500000,-3.500000);
\draw[-triangle 60,thick](0.500000,-5.500000) -- (5.500000,-4.500000);
\draw[-triangle 60,thick](0.500000,-6.500000) -- (2.500000,-5.500000);
\draw[-triangle 60,thick](0.500000,-7.500000) -- (1.500000,-6.500000);
\end{tikzpicture}
}\scalebox{0.380000}{\begin{tikzpicture}[turtle/distance=1cm]\filldraw [thick,white,fill=gray!25,turtle={home,right,forward,left,right,right,left,forward=8,right,forward,forward,right,forward=6,left,forward,right,forward=1,left,forward,left,forward=4,right,forward,right,forward=3,left,forward,right,forward=1,left,forward,left,forward=3,right,forward,right,forward=4,right,forward=8}];\seeddos{1.000000}{0}\nonseeddos{2.000000}{0}\seeddos{3.000000}{0}\nonseeddos{4.000000}{0}\nonseeddos{5.000000}{0}\nonseeddos{6.000000}{0}\nonseeddos{7.000000}{0}\transparentnonseed{8}{0}\seeddos{1.000000}{1}\nonseeddos{2.000000}{1}\seeddos{3.000000}{1}\nonseeddos{4.000000}{1}\nonseeddos{5.000000}{1}\nonseeddos{6.000000}{1}\nonseeddos{7.000000}{1}\transparentnonseed{8}{1}\seeddos{1.000000}{2}\transparentnonseed{2}{2}\transparentseed{1}{3}\nonseeddos{1.000000}{4}\nonseeddos{2.000000}{4}\nonseeddos{3.000000}{4}\nonseeddos{4.000000}{4}\transparentnonseed{5}{4}\seeddos{1.000000}{5}\transparentseed{2}{5}\transparentseed{1}{6}\seeddos{1.000000}{7}\seeddos{2.000000}{7}\end{tikzpicture}
}\end{center}Set the last two elments in the last row of the table as old-order new seeds.

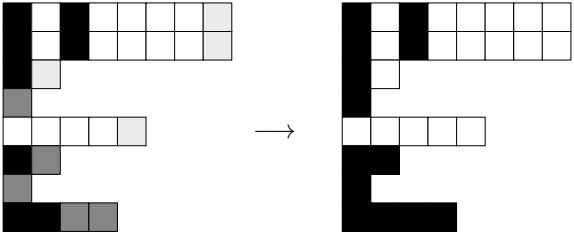
\begin{center}\scalebox{0.380000}{\begin{tikzpicture}[turtle/distance=1cm]\filldraw [thick,white,fill=gray!25,turtle={home,right,forward,left,right,right,left,forward=8,right,forward,forward,right,forward=6,left,forward,right,forward=1,left,forward,left,forward=4,right,forward,right,forward=3,left,forward,right,forward=1,left,forward,left,forward=3,right,forward,right,forward=4,right,forward=8}];\netejaseed{1.000000}{0}\netejanonseed{2.000000}{0}\netejaseed{3.000000}{0}\netejanonseed{4.000000}{0}\netejanonseed{5.000000}{0}\netejanonseed{6.000000}{0}\netejanonseed{7.000000}{0}\netejatransparentnonseed{8}{0}\netejaseed{1.000000}{1}\netejanonseed{2.000000}{1}\netejaseed{3.000000}{1}\netejanonseed{4.000000}{1}\netejanonseed{5.000000}{1}\netejanonseed{6.000000}{1}\netejanonseed{7.000000}{1}\netejatransparentnonseed{8}{1}\netejaseed{1.000000}{2}\netejatransparentnonseed{2}{2}\netejatransparentseed{1}{3}\netejanonseed{1.000000}{4}\netejanonseed{2.000000}{4}\netejanonseed{3.000000}{4}\netejanonseed{4.000000}{4}\netejatransparentnonseed{5}{4}\netejaseed{1.000000}{5}\netejatransparentseed{2}{5}\netejatransparentseed{1}{6}\netejaseed{1.000000}{7}\netejaseed{2.000000}{7}\netejatransparentseed{3}{7}\netejatransparentseed{4}{7}\end{tikzpicture}
}\ \ \raisebox{1.25cm}{$\longrightarrow$}\ \ \scalebox{0.380000}{\begin{tikzpicture}[turtle/distance=1cm]\filldraw [thick,white,fill=gray!25,turtle={home,right,forward,left,right,right,left,forward=8,right,forward,forward,right,forward=6,left,forward,right,forward=1,left,forward,left,forward=4,right,forward,right,forward=3,left,forward,right,forward=1,left,forward,left,forward=3,right,forward,right,forward=4,right,forward=8}];\seeddos{1.000000}{0}\nonseeddos{2.000000}{0}\seeddos{3.000000}{0}\nonseeddos{4.000000}{0}\nonseeddos{5.000000}{0}\nonseeddos{6.000000}{0}\nonseeddos{7.000000}{0}\nonseeddos{8.000000}{0}\seeddos{1.000000}{1}\nonseeddos{2.000000}{1}\seeddos{3.000000}{1}\nonseeddos{4.000000}{1}\nonseeddos{5.000000}{1}\nonseeddos{6.000000}{1}\nonseeddos{7.000000}{1}\nonseeddos{8.000000}{1}\seeddos{1.000000}{2}\nonseeddos{2.000000}{2}\seeddos{1.000000}{3}\nonseeddos{1.000000}{4}\nonseeddos{2.000000}{4}\nonseeddos{3.000000}{4}\nonseeddos{4.000000}{4}\nonseeddos{5.000000}{4}\seeddos{1.000000}{5}\seeddos{2.000000}{5}\seeddos{1.000000}{6}\seeddos{1.000000}{7}\seeddos{2.000000}{7}\seeddos{3.000000}{7}\seeddos{4.000000}{7}\end{tikzpicture}
}\end{center}Suppose that now we want to take away the generator c+1=30+1=31.\begin{center}\scalebox{0.380000}
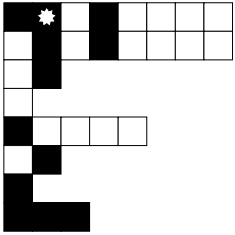
{\begin{tikzpicture}[turtle/distance=1cm]\filldraw [thick,white,fill=gray!25,turtle={home,right,forward,left,right,right,left,forward=8,right,forward,forward,right,forward=6,left,forward,right,forward=1,left,forward,left,forward=4,right,forward,right,forward=3,left,forward,right,forward=1,left,forward,left,forward=2,right,forward,right,forward=3,right,forward=8}];\seeddos{1.000000}{0}\seeddos{2.000000}{0}\nonseeddos{3.000000}{0}\seeddos{4.000000}{0}\nonseeddos{5.000000}{0}\nonseeddos{6.000000}{0}\nonseeddos{7.000000}{0}\nonseeddos{8.000000}{0}\nonseeddos{1.000000}{1}\seeddos{2.000000}{1}\nonseeddos{3.000000}{1}\seeddos{4.000000}{1}\nonseeddos{5.000000}{1}\nonseeddos{6.000000}{1}\nonseeddos{7.000000}{1}\nonseeddos{8.000000}{1}\nonseeddos{1.000000}{2}\seeddos{2.000000}{2}\nonseeddos{1.000000}{3}\seeddos{1.000000}{4}\nonseeddos{2.000000}{4}\nonseeddos{3.000000}{4}\nonseeddos{4.000000}{4}\nonseeddos{5.000000}{4}\nonseeddos{1.000000}{5}\seeddos{2.000000}{5}\seeddos{1.000000}{6}\seeddos{1.000000}{7}\seeddos{2.000000}{7}\seeddos{3.000000}{7}\node[draw, fill=white, star, star points=9, minimum size=7mm] at (2.500000,-.5){};
\end{tikzpicture}
}\end{center}Draw the contour of the new table of seeds.

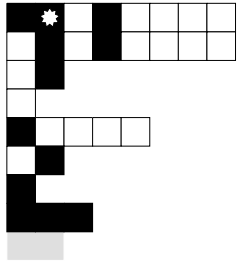
\begin{center}\scalebox{0.380000}{\begin{tikzpicture}[turtle/distance=1cm]\filldraw [thick,white,fill=gray!25,turtle={home,right,forward,left,right,right,left,forward=8,right,forward,forward,right,forward=6,left,forward,right,forward=1,left,forward,left,forward=4,right,forward,right,forward=3,left,forward,right,forward=1,left,forward,left,forward=2,right,forward,right,forward=1,left,forward,right,forward=2,right,forward=9}];\seeddos{1.000000}{0}\seeddos{2.000000}{0}\nonseeddos{3.000000}{0}\seeddos{4.000000}{0}\nonseeddos{5.000000}{0}\nonseeddos{6.000000}{0}\nonseeddos{7.000000}{0}\nonseeddos{8.000000}{0}\nonseeddos{1.000000}{1}\seeddos{2.000000}{1}\nonseeddos{3.000000}{1}\seeddos{4.000000}{1}\nonseeddos{5.000000}{1}\nonseeddos{6.000000}{1}\nonseeddos{7.000000}{1}\nonseeddos{8.000000}{1}\nonseeddos{1.000000}{2}\seeddos{2.000000}{2}\nonseeddos{1.000000}{3}\seeddos{1.000000}{4}\nonseeddos{2.000000}{4}\nonseeddos{3.000000}{4}\nonseeddos{4.000000}{4}\nonseeddos{5.000000}{4}\nonseeddos{1.000000}{5}\seeddos{2.000000}{5}\seeddos{1.000000}{6}\seeddos{1.000000}{7}\seeddos{2.000000}{7}\seeddos{3.000000}{7}\node[draw, fill=white, star, star points=9, minimum size=7mm] at (2.500000,-.5){};
\end{tikzpicture}
}\end{center}Discard the values corresponding to elements that are smaller than the new Frobenius number, keep shadowed the values corresponding to the new Frobenius number.

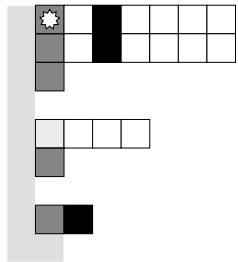
\begin{center}\scalebox{0.380000}{\begin{tikzpicture}[turtle/distance=1cm]\filldraw [thick,white,fill=gray!25,turtle={home,right,forward,left,right,right,left,forward=8,right,forward,forward,right,forward=6,left,forward,right,forward=1,left,forward,left,forward=4,right,forward,right,forward=3,left,forward,right,forward=1,left,forward,left,forward=2,right,forward,right,forward=1,left,forward,right,forward=2,right,forward=9}];\transparentseed{2.000000}{0}\nonseeddos{3.000000}{0}\seeddos{4.000000}{0}\nonseeddos{5.000000}{0}\nonseeddos{6.000000}{0}\nonseeddos{7.000000}{0}\nonseeddos{8.000000}{0}\transparentseed{2.000000}{1}\nonseeddos{3.000000}{1}\seeddos{4.000000}{1}\nonseeddos{5.000000}{1}\nonseeddos{6.000000}{1}\nonseeddos{7.000000}{1}\nonseeddos{8.000000}{1}\transparentseed{2.000000}{2}\transparentnonseed{2.000000}{4}\nonseeddos{3.000000}{4}\nonseeddos{4.000000}{4}\nonseeddos{5.000000}{4}\transparentseed{2.000000}{5}\transparentseed{2.000000}{7}\seeddos{3.000000}{7}\node[draw, fill=white, star, star points=9, minimum size=7mm] at (2.500000,-.5){};
\end{tikzpicture}
}\end{center}Move the old values to the left of the table and fix the old-order recycled seeds.

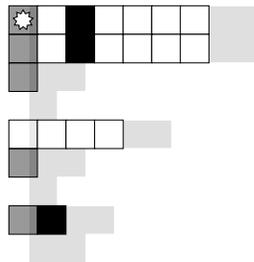
\begin{center}\scalebox{0.380000}{\begin{tikzpicture}[turtle/distance=1cm]\filldraw [thick,white,fill=gray!25,turtle={home,right,forward,left,right,right,left,forward=8,right,forward,forward,right,forward=6,left,forward,right,forward=1,left,forward,left,forward=4,right,forward,right,forward=3,left,forward,right,forward=1,left,forward,left,forward=2,right,forward,right,forward=1,left,forward,right,forward=2,right,forward=9}];\transparentseed{0.300000}{0}\nonseeddos{1.300000}{0}\seeddos{2.300000}{0}\nonseeddos{3.300000}{0}\nonseeddos{4.300000}{0}\nonseeddos{5.300000}{0}\nonseeddos{6.300000}{0}\transparentseed{0.300000}{1}\nonseeddos{1.300000}{1}\seeddos{2.300000}{1}\nonseeddos{3.300000}{1}\nonseeddos{4.300000}{1}\nonseeddos{5.300000}{1}\nonseeddos{6.300000}{1}\transparentseed{0.300000}{2}\transparentnonseed{0.300000}{4}\nonseeddos{1.300000}{4}\nonseeddos{2.300000}{4}\nonseeddos{3.300000}{4}\transparentseed{0.300000}{5}\transparentseed{0.300000}{7}\seeddos{1.300000}{7}\node[draw, fill=white, star, star points=9, minimum size=7mm] at (0.800000,-.5){};
\end{tikzpicture}
}\end{center}Obtain the old-order new seeds.

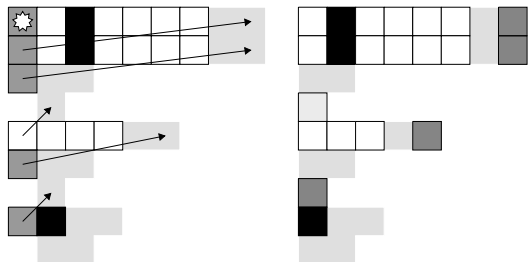
\begin{center}\scalebox{0.380000}{\begin{tikzpicture}[turtle/distance=1cm]\filldraw [thick,white,fill=gray!25,turtle={home,right,forward,left,right,right,left,forward=8,right,forward,forward,right,forward=6,left,forward,right,forward=1,left,forward,left,forward=4,right,forward,right,forward=3,left,forward,right,forward=1,left,forward,left,forward=2,right,forward,right,forward=1,left,forward,right,forward=2,right,forward=9}];\transparentseed{0.000000}{0}\nonseeddos{1.000000}{0}\seeddos{2.000000}{0}\nonseeddos{3.000000}{0}\nonseeddos{4.000000}{0}\nonseeddos{5.000000}{0}\nonseeddos{6.000000}{0}\transparentseed{0.000000}{1}\nonseeddos{1.000000}{1}\seeddos{2.000000}{1}\nonseeddos{3.000000}{1}\nonseeddos{4.000000}{1}\nonseeddos{5.000000}{1}\nonseeddos{6.000000}{1}\transparentseed{0.000000}{2}\transparentnonseed{0.000000}{4}\nonseeddos{1.000000}{4}\nonseeddos{2.000000}{4}\nonseeddos{3.000000}{4}\transparentseed{0.000000}{5}\transparentseed{0.000000}{7}\seeddos{1.000000}{7}\node[draw, fill=white, star, star points=9, minimum size=7mm] at (0.500000,-.5){};
\draw[-triangle 60,thick](0.500000,-1.500000) -- (8.500000,-0.500000);
\draw[-triangle 60,thick](0.500000,-2.500000) -- (8.500000,-1.500000);
\draw[-triangle 60,thick](0.500000,-4.500000) -- (1.500000,-3.500000);
\draw[-triangle 60,thick](0.500000,-5.500000) -- (5.500000,-4.500000);
\draw[-triangle 60,thick](0.500000,-7.500000) -- (1.500000,-6.500000);
\end{tikzpicture}
}\scalebox{0.380000}{\begin{tikzpicture}[turtle/distance=1cm]\filldraw [thick,white,fill=gray!25,turtle={home,right,forward,left,right,right,left,forward=8,right,forward,forward,right,forward=6,left,forward,right,forward=1,left,forward,left,forward=4,right,forward,right,forward=3,left,forward,right,forward=1,left,forward,left,forward=2,right,forward,right,forward=1,left,forward,right,forward=2,right,forward=9}];\nonseeddos{1.000000}{0}\seeddos{2.000000}{0}\nonseeddos{3.000000}{0}\nonseeddos{4.000000}{0}\nonseeddos{5.000000}{0}\nonseeddos{6.000000}{0}\transparentseed{8}{0}\nonseeddos{1.000000}{1}\seeddos{2.000000}{1}\nonseeddos{3.000000}{1}\nonseeddos{4.000000}{1}\nonseeddos{5.000000}{1}\nonseeddos{6.000000}{1}\transparentseed{8}{1}\transparentnonseed{1}{3}\nonseeddos{1.000000}{4}\nonseeddos{2.000000}{4}\nonseeddos{3.000000}{4}\transparentseed{5}{4}\transparentseed{1}{6}\seeddos{1.000000}{7}\end{tikzpicture}
}\end{center}Set the last elment in the last but one row of the table as one old-order new seed and set the two elements in the last row as two new-order seeds.

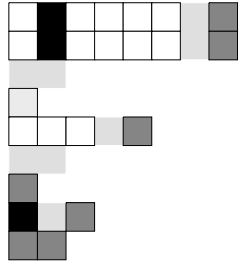
\begin{center}\scalebox{0.380000}{\begin{tikzpicture}[turtle/distance=1cm]\filldraw [thick,white,fill=gray!25,turtle={home,right,forward,left,right,right,left,forward=8,right,forward,forward,right,forward=6,left,forward,right,forward=1,left,forward,left,forward=4,right,forward,right,forward=3,left,forward,right,forward=1,left,forward,left,forward=2,right,forward,right,forward=1,left,forward,right,forward=2,right,forward=9}];\netejanonseed{1.000000}{0}\netejaseed{2.000000}{0}\netejanonseed{3.000000}{0}\netejanonseed{4.000000}{0}\netejanonseed{5.000000}{0}\netejanonseed{6.000000}{0}\netejatransparentseed{8}{0}\netejanonseed{1.000000}{1}\netejaseed{2.000000}{1}\netejanonseed{3.000000}{1}\netejanonseed{4.000000}{1}\netejanonseed{5.000000}{1}\netejanonseed{6.000000}{1}\netejatransparentseed{8}{1}\netejatransparentnonseed{1}{3}\netejanonseed{1.000000}{4}\netejanonseed{2.000000}{4}\netejanonseed{3.000000}{4}\netejatransparentseed{5}{4}\netejatransparentseed{1}{6}\netejaseed{1.000000}{7}\netejatransparentseed{3}{7}\netejatransparentseed{1}{8}\netejatransparentseed{2}{8}\end{tikzpicture}
}\end{center}The remaining empty boxes are non-seeds.

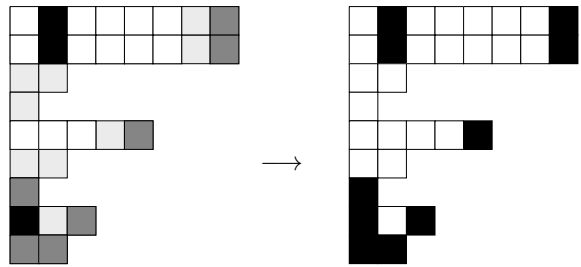
\begin{center}\scalebox{0.380000}{\begin{tikzpicture}[turtle/distance=1cm]\filldraw [thick,white,fill=gray!25,turtle={home,right,forward,left,right,right,left,forward=8,right,forward,forward,right,forward=6,left,forward,right,forward=1,left,forward,left,forward=4,right,forward,right,forward=3,left,forward,right,forward=1,left,forward,left,forward=2,right,forward,right,forward=1,left,forward,right,forward=2,right,forward=9}];\transparentnonseed{1}{0}\transparentnonseed{2}{0}\transparentnonseed{3}{0}\transparentnonseed{4}{0}\transparentnonseed{5}{0}\transparentnonseed{6}{0}\transparentnonseed{7}{0}\transparentnonseed{8}{0}\transparentnonseed{1}{1}\transparentnonseed{2}{1}\transparentnonseed{3}{1}\transparentnonseed{4}{1}\transparentnonseed{5}{1}\transparentnonseed{6}{1}\transparentnonseed{7}{1}\transparentnonseed{8}{1}\transparentnonseed{1}{2}\transparentnonseed{2}{2}\transparentnonseed{1}{3}\transparentnonseed{1}{4}\transparentnonseed{2}{4}\transparentnonseed{3}{4}\transparentnonseed{4}{4}\transparentnonseed{5}{4}\transparentnonseed{1}{5}\transparentnonseed{2}{5}\transparentnonseed{1}{6}\transparentnonseed{1}{7}\transparentnonseed{2}{7}\transparentnonseed{3}{7}\transparentnonseed{1}{8}\transparentnonseed{2}{8}\netejanonseed{1.000000}{0}\netejaseed{2.000000}{0}\netejanonseed{3.000000}{0}\netejanonseed{4.000000}{0}\netejanonseed{5.000000}{0}\netejanonseed{6.000000}{0}\netejatransparentseed{8}{0}\netejanonseed{1.000000}{1}\netejaseed{2.000000}{1}\netejanonseed{3.000000}{1}\netejanonseed{4.000000}{1}\netejanonseed{5.000000}{1}\netejanonseed{6.000000}{1}\netejatransparentseed{8}{1}\netejatransparentnonseed{1}{3}\netejanonseed{1.000000}{4}\netejanonseed{2.000000}{4}\netejanonseed{3.000000}{4}\netejatransparentseed{5}{4}\netejatransparentseed{1}{6}\netejaseed{1.000000}{7}\netejatransparentseed{3}{7}\netejatransparentseed{1}{8}\netejatransparentseed{2}{8}\end{tikzpicture}
}\ \ \raisebox{1.25cm}{$\longrightarrow$}\ \ \scalebox{0.380000}{\begin{tikzpicture}[turtle/distance=1cm]\filldraw [thick,white,fill=gray!25,turtle={home,right,forward,left,right,right,left,forward=8,right,forward,forward,right,forward=6,left,forward,right,forward=1,left,forward,left,forward=4,right,forward,right,forward=3,left,forward,right,forward=1,left,forward,left,forward=2,right,forward,right,forward=1,left,forward,right,forward=2,right,forward=9}];\nonseeddos{1.000000}{0}\seeddos{2.000000}{0}\nonseeddos{3.000000}{0}\nonseeddos{4.000000}{0}\nonseeddos{5.000000}{0}\nonseeddos{6.000000}{0}\nonseeddos{7.000000}{0}\seeddos{8.000000}{0}\nonseeddos{1.000000}{1}\seeddos{2.000000}{1}\nonseeddos{3.000000}{1}\nonseeddos{4.000000}{1}\nonseeddos{5.000000}{1}\nonseeddos{6.000000}{1}\nonseeddos{7.000000}{1}\seeddos{8.000000}{1}\nonseeddos{1.000000}{2}\nonseeddos{2.000000}{2}\nonseeddos{1.000000}{3}\nonseeddos{1.000000}{4}\nonseeddos{2.000000}{4}\nonseeddos{3.000000}{4}\nonseeddos{4.000000}{4}\seeddos{5.000000}{4}\nonseeddos{1.000000}{5}\nonseeddos{2.000000}{5}\seeddos{1.000000}{6}\seeddos{1.000000}{7}\nonseeddos{2.000000}{7}\seeddos{3.000000}{7}\seeddos{1.000000}{8}\seeddos{2.000000}{8}\end{tikzpicture}
}\end{center}Suppose that now we want to take away the generator c+3=30+3=33.\begin{center}\scalebox{0.380000}
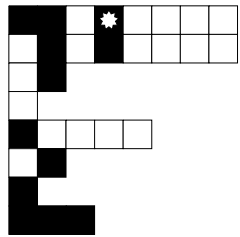
{\begin{tikzpicture}[turtle/distance=1cm]\filldraw [thick,white,fill=gray!25,turtle={home,right,forward,left,right,right,left,forward=8,right,forward,forward,right,forward=6,left,forward,right,forward=1,left,forward,left,forward=4,right,forward,right,forward=3,left,forward,right,forward=1,left,forward,left,forward=2,right,forward,right,forward=3,right,forward=8}];\seeddos{1.000000}{0}\seeddos{2.000000}{0}\nonseeddos{3.000000}{0}\seeddos{4.000000}{0}\nonseeddos{5.000000}{0}\nonseeddos{6.000000}{0}\nonseeddos{7.000000}{0}\nonseeddos{8.000000}{0}\nonseeddos{1.000000}{1}\seeddos{2.000000}{1}\nonseeddos{3.000000}{1}\seeddos{4.000000}{1}\nonseeddos{5.000000}{1}\nonseeddos{6.000000}{1}\nonseeddos{7.000000}{1}\nonseeddos{8.000000}{1}\nonseeddos{1.000000}{2}\seeddos{2.000000}{2}\nonseeddos{1.000000}{3}\seeddos{1.000000}{4}\nonseeddos{2.000000}{4}\nonseeddos{3.000000}{4}\nonseeddos{4.000000}{4}\nonseeddos{5.000000}{4}\nonseeddos{1.000000}{5}\seeddos{2.000000}{5}\seeddos{1.000000}{6}\seeddos{1.000000}{7}\seeddos{2.000000}{7}\seeddos{3.000000}{7}\node[draw, fill=white, star, star points=9, minimum size=7mm] at (4.500000,-.5){};
\end{tikzpicture}
}\end{center}Draw the contour of the new table of seeds.

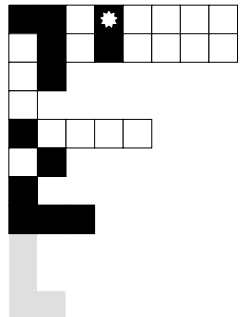
\begin{center}\scalebox{0.380000}{\begin{tikzpicture}[turtle/distance=1cm]\filldraw [thick,white,fill=gray!25,turtle={home,right,forward,left,right,right,left,forward=8,right,forward,forward,right,forward=6,left,forward,right,forward=1,left,forward,left,forward=4,right,forward,right,forward=3,left,forward,right,forward=1,left,forward,left,forward=2,right,forward,right,forward=2,left,forward,forward,left,forward=1,right,forward,right,forward=2,right,forward=11}];\seeddos{1.000000}{0}\seeddos{2.000000}{0}\nonseeddos{3.000000}{0}\seeddos{4.000000}{0}\nonseeddos{5.000000}{0}\nonseeddos{6.000000}{0}\nonseeddos{7.000000}{0}\nonseeddos{8.000000}{0}\nonseeddos{1.000000}{1}\seeddos{2.000000}{1}\nonseeddos{3.000000}{1}\seeddos{4.000000}{1}\nonseeddos{5.000000}{1}\nonseeddos{6.000000}{1}\nonseeddos{7.000000}{1}\nonseeddos{8.000000}{1}\nonseeddos{1.000000}{2}\seeddos{2.000000}{2}\nonseeddos{1.000000}{3}\seeddos{1.000000}{4}\nonseeddos{2.000000}{4}\nonseeddos{3.000000}{4}\nonseeddos{4.000000}{4}\nonseeddos{5.000000}{4}\nonseeddos{1.000000}{5}\seeddos{2.000000}{5}\seeddos{1.000000}{6}\seeddos{1.000000}{7}\seeddos{2.000000}{7}\seeddos{3.000000}{7}\node[draw, fill=white, star, star points=9, minimum size=7mm] at (4.500000,-.5){};
\end{tikzpicture}
}\end{center}Discard the values corresponding to elements that are smaller than the new Frobenius number, keep shadowed the values corresponding to the new Frobenius number.

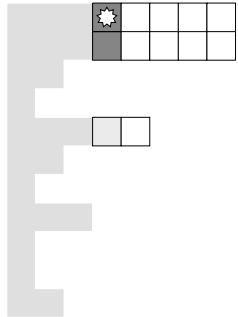
\begin{center}\scalebox{0.380000}{\begin{tikzpicture}[turtle/distance=1cm]\filldraw [thick,white,fill=gray!25,turtle={home,right,forward,left,right,right,left,forward=8,right,forward,forward,right,forward=6,left,forward,right,forward=1,left,forward,left,forward=4,right,forward,right,forward=3,left,forward,right,forward=1,left,forward,left,forward=2,right,forward,right,forward=2,left,forward,forward,left,forward=1,right,forward,right,forward=2,right,forward=11}];\transparentseed{4.000000}{0}\nonseeddos{5.000000}{0}\nonseeddos{6.000000}{0}\nonseeddos{7.000000}{0}\nonseeddos{8.000000}{0}\transparentseed{4.000000}{1}\nonseeddos{5.000000}{1}\nonseeddos{6.000000}{1}\nonseeddos{7.000000}{1}\nonseeddos{8.000000}{1}\transparentnonseed{4.000000}{4}\nonseeddos{5.000000}{4}\node[draw, fill=white, star, star points=9, minimum size=7mm] at (4.500000,-.5){};
\end{tikzpicture}
}\end{center}Move the old values to the left of the table and fix the old-order recycled seeds.

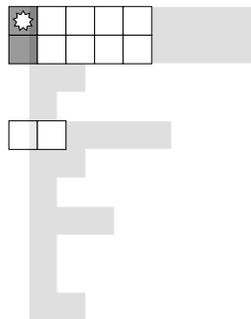
\begin{center}\scalebox{0.380000}{\begin{tikzpicture}[turtle/distance=1cm]\filldraw [thick,white,fill=gray!25,turtle={home,right,forward,left,right,right,left,forward=8,right,forward,forward,right,forward=6,left,forward,right,forward=1,left,forward,left,forward=4,right,forward,right,forward=3,left,forward,right,forward=1,left,forward,left,forward=2,right,forward,right,forward=2,left,forward,forward,left,forward=1,right,forward,right,forward=2,right,forward=11}];\transparentseed{0.300000}{0}\nonseeddos{1.300000}{0}\nonseeddos{2.300000}{0}\nonseeddos{3.300000}{0}\nonseeddos{4.300000}{0}\transparentseed{0.300000}{1}\nonseeddos{1.300000}{1}\nonseeddos{2.300000}{1}\nonseeddos{3.300000}{1}\nonseeddos{4.300000}{1}\transparentnonseed{0.300000}{4}\nonseeddos{1.300000}{4}\node[draw, fill=white, star, star points=9, minimum size=7mm] at (0.800000,-.5){};
\end{tikzpicture}
}\end{center}Obtain the old-order new seeds.

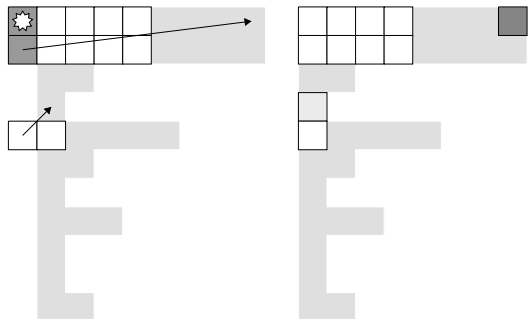
\begin{center}\scalebox{0.380000}{\begin{tikzpicture}[turtle/distance=1cm]\filldraw [thick,white,fill=gray!25,turtle={home,right,forward,left,right,right,left,forward=8,right,forward,forward,right,forward=6,left,forward,right,forward=1,left,forward,left,forward=4,right,forward,right,forward=3,left,forward,right,forward=1,left,forward,left,forward=2,right,forward,right,forward=2,left,forward,forward,left,forward=1,right,forward,right,forward=2,right,forward=11}];\transparentseed{0.000000}{0}\nonseeddos{1.000000}{0}\nonseeddos{2.000000}{0}\nonseeddos{3.000000}{0}\nonseeddos{4.000000}{0}\transparentseed{0.000000}{1}\nonseeddos{1.000000}{1}\nonseeddos{2.000000}{1}\nonseeddos{3.000000}{1}\nonseeddos{4.000000}{1}\transparentnonseed{0.000000}{4}\nonseeddos{1.000000}{4}\node[draw, fill=white, star, star points=9, minimum size=7mm] at (0.500000,-.5){};
\draw[-triangle 60,thick](0.500000,-1.500000) -- (8.500000,-0.500000);
\draw[-triangle 60,thick](0.500000,-4.500000) -- (1.500000,-3.500000);
\end{tikzpicture}
}\scalebox{0.380000}{\begin{tikzpicture}[turtle/distance=1cm]\filldraw [thick,white,fill=gray!25,turtle={home,right,forward,left,right,right,left,forward=8,right,forward,forward,right,forward=6,left,forward,right,forward=1,left,forward,left,forward=4,right,forward,right,forward=3,left,forward,right,forward=1,left,forward,left,forward=2,right,forward,right,forward=2,left,forward,forward,left,forward=1,right,forward,right,forward=2,right,forward=11}];\nonseeddos{1.000000}{0}\nonseeddos{2.000000}{0}\nonseeddos{3.000000}{0}\nonseeddos{4.000000}{0}\transparentseed{8}{0}\nonseeddos{1.000000}{1}\nonseeddos{2.000000}{1}\nonseeddos{3.000000}{1}\nonseeddos{4.000000}{1}\transparentnonseed{1}{3}\nonseeddos{1.000000}{4}\end{tikzpicture}
}\end{center}Set the unique element in the last but one row of the table and the two elements in the last row of the table as new-order seeds.

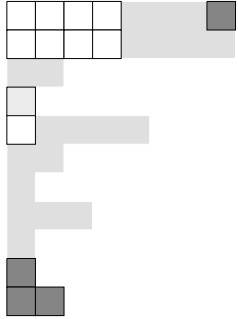
\begin{center}\scalebox{0.380000}{\begin{tikzpicture}[turtle/distance=1cm]\filldraw [thick,white,fill=gray!25,turtle={home,right,forward,left,right,right,left,forward=8,right,forward,forward,right,forward=6,left,forward,right,forward=1,left,forward,left,forward=4,right,forward,right,forward=3,left,forward,right,forward=1,left,forward,left,forward=2,right,forward,right,forward=2,left,forward,forward,left,forward=1,right,forward,right,forward=2,right,forward=11}];\netejanonseed{1.000000}{0}\netejanonseed{2.000000}{0}\netejanonseed{3.000000}{0}\netejanonseed{4.000000}{0}\netejatransparentseed{8}{0}\netejanonseed{1.000000}{1}\netejanonseed{2.000000}{1}\netejanonseed{3.000000}{1}\netejanonseed{4.000000}{1}\netejatransparentnonseed{1}{3}\netejanonseed{1.000000}{4}\netejatransparentseed{1}{9}\netejatransparentseed{1}{10}\netejatransparentseed{2}{10}\end{tikzpicture}
}\end{center}The remaining empty boxes are non-seeds.

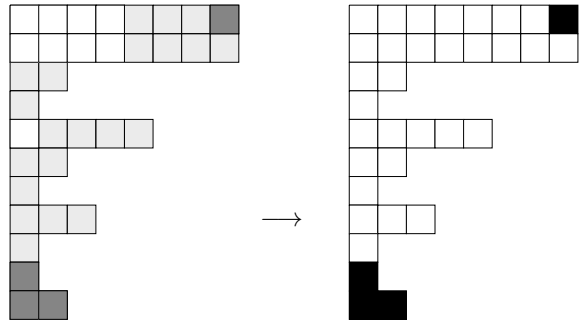
\begin{center}\scalebox{0.380000}{\begin{tikzpicture}[turtle/distance=1cm]\filldraw [thick,white,fill=gray!25,turtle={home,right,forward,left,right,right,left,forward=8,right,forward,forward,right,forward=6,left,forward,right,forward=1,left,forward,left,forward=4,right,forward,right,forward=3,left,forward,right,forward=1,left,forward,left,forward=2,right,forward,right,forward=2,left,forward,forward,left,forward=1,right,forward,right,forward=2,right,forward=11}];\transparentnonseed{1}{0}\transparentnonseed{2}{0}\transparentnonseed{3}{0}\transparentnonseed{4}{0}\transparentnonseed{5}{0}\transparentnonseed{6}{0}\transparentnonseed{7}{0}\transparentnonseed{8}{0}\transparentnonseed{1}{1}\transparentnonseed{2}{1}\transparentnonseed{3}{1}\transparentnonseed{4}{1}\transparentnonseed{5}{1}\transparentnonseed{6}{1}\transparentnonseed{7}{1}\transparentnonseed{8}{1}\transparentnonseed{1}{2}\transparentnonseed{2}{2}\transparentnonseed{1}{3}\transparentnonseed{1}{4}\transparentnonseed{2}{4}\transparentnonseed{3}{4}\transparentnonseed{4}{4}\transparentnonseed{5}{4}\transparentnonseed{1}{5}\transparentnonseed{2}{5}\transparentnonseed{1}{6}\transparentnonseed{1}{7}\transparentnonseed{2}{7}\transparentnonseed{3}{7}\transparentnonseed{1}{8}\transparentnonseed{1}{9}\transparentnonseed{1}{10}\transparentnonseed{2}{10}\netejanonseed{1.000000}{0}\netejanonseed{2.000000}{0}\netejanonseed{3.000000}{0}\netejanonseed{4.000000}{0}\netejatransparentseed{8}{0}\netejanonseed{1.000000}{1}\netejanonseed{2.000000}{1}\netejanonseed{3.000000}{1}\netejanonseed{4.000000}{1}\netejatransparentnonseed{1}{3}\netejanonseed{1.000000}{4}\netejatransparentseed{1}{9}\netejatransparentseed{1}{10}\netejatransparentseed{2}{10}\end{tikzpicture}
}\ \ \raisebox{1.25cm}{$\longrightarrow$}\ \ \scalebox{0.380000}{\begin{tikzpicture}[turtle/distance=1cm]\filldraw [thick,white,fill=gray!25,turtle={home,right,forward,left,right,right,left,forward=8,right,forward,forward,right,forward=6,left,forward,right,forward=1,left,forward,left,forward=4,right,forward,right,forward=3,left,forward,right,forward=1,left,forward,left,forward=2,right,forward,right,forward=2,left,forward,forward,left,forward=1,right,forward,right,forward=2,right,forward=11}];\nonseeddos{1.000000}{0}\nonseeddos{2.000000}{0}\nonseeddos{3.000000}{0}\nonseeddos{4.000000}{0}\nonseeddos{5.000000}{0}\nonseeddos{6.000000}{0}\nonseeddos{7.000000}{0}\seeddos{8.000000}{0}\nonseeddos{1.000000}{1}\nonseeddos{2.000000}{1}\nonseeddos{3.000000}{1}\nonseeddos{4.000000}{1}\nonseeddos{5.000000}{1}\nonseeddos{6.000000}{1}\nonseeddos{7.000000}{1}\nonseeddos{8.000000}{1}\nonseeddos{1.000000}{2}\nonseeddos{2.000000}{2}\nonseeddos{1.000000}{3}\nonseeddos{1.000000}{4}\nonseeddos{2.000000}{4}\nonseeddos{3.000000}{4}\nonseeddos{4.000000}{4}\nonseeddos{5.000000}{4}\nonseeddos{1.000000}{5}\nonseeddos{2.000000}{5}\nonseeddos{1.000000}{6}\nonseeddos{1.000000}{7}\nonseeddos{2.000000}{7}\nonseeddos{3.000000}{7}\nonseeddos{1.000000}{8}\seeddos{1.000000}{9}\seeddos{1.000000}{10}\seeddos{2.000000}{10}\end{tikzpicture}
}\end{center}

\section{Eliahou semigroups and Wilf conjecture verification extended up to genus 65}
\label{s:EliahouWilf}

Fix a numerical semigroup $\Lambda$ with conductor $c$ and multiplicity $m$.
Let $q=\lceil\frac{c}{m}\rceil$, and let $\rho=qm-c$ be the remainder of the division of $c$ by $m$.
Suppose that the left elements of $\Lambda$ are $\{\lambda_0,\dots,\lambda_L\}$ and suppose that $P$ is the set of primitive elements of $\Lambda$.
The Eliahou constant is defined as
$$E(\Lambda)=\#(P\cap \{\lambda_0,\dots,\lambda_L\})L-q(m-\#(P\setminus \{\lambda_0,\dots,\lambda_L\}))+\rho.$$
Shalom Eliahou proved that if $E(\Lambda)\geq 0$ then $\Lambda$ satifies the Wilf conjecture \cite{eliahou}.
Semigroups for which the Eliahou constant is negative are very unusual. We will denote them Eliahou semigroups.
According to the same reference, it was computed by Jean Fromentin that the unique Eliahou semigroups of genus $g\leq 60$ are exactly
\begin{itemize}
  \item
    $\langle 14,22,23\rangle\mid_{56}$,

    \bigskip

  \item    $\langle 16,25,26\rangle\mid_{64}$,

    \bigskip

  \item  $\langle 17,26,28\rangle\mid_{68}$,

    \bigskip

  \item
    
    $\langle 17,27,28\rangle\mid_{68}$,

    \bigskip

    \item
      $\langle 18,28,29\rangle\mid_{72}$
\end{itemize}
for genus 43, 51, 55, 55 and 59, respectively,
where $\langle a,b,c\rangle\mid_{\kappa}$ means the minimum semigroup containing $a,b,c$ and all integers larger than or equal to $\kappa$.

Using a parallelized version of the seeds algorithm we found that the unique Eliahou semigroups with genus between $61$ and $65$ are exactly,

$$\{0,19, 29, 31, 38, 48, 50, 57, 58, 60, 62, 67, 69, 76,\dots\}=\langle 19,29,31\rangle\mid_{76}$$
and
$$\{0,19, 30, 31, 38, 49, 50, 57, 60, 61, 62, 68, 69, 76,\dots\}=\langle 19,30,31\rangle\mid_{76}$$

Using $p$ for the number of primitive elements and $r$ for the number of right primitive elements, the parameters of these semigroups are

\begin{itemize}
%\item
%  $\{0,14, 22, 23, 28, 36, 37, 42, 44, 45, 46, 50, 51, 56,\dots\}$
%\\($g=43, c=56, m=14, p=7, r=4, L=13, q=4, \rho=0$)
%
%\bigskip
%
%\item
%$\{0,16, 25, 26, 32, 41, 42, 48, 50, 51, 52, 57, 58, 64,\dots\}$
%\\($g=51, c=64, m=16, p=9, r=6, L=13, q=4, \rho=0$)
%
%\bigskip
%
%\item
%$\{0,17, 26, 28, 34, 43, 45, 51, 52, 54, 56, 60, 62, 68,\dots\}$
%\\($g=55, c=68, m=17, p=10, r=7, L=13, q=4, \rho=0$)
%
%\bigskip
%
%\item 
%$\{0,17, 27, 28, 34, 44, 45, 51, 54, 55, 56, 61, 62, 68,\dots\}$
%\\($g=55, c=68, m=17, p=10, r=7, L=13, q=4, \rho=0$)
%
%\bigskip
%
%\item 
%$\{0,  18, 28, 29, 36, 46, 47, 54, 56, 57, 58, 64, 65, 72,\dots\}$
%\\($g=59, c=72, m=18, p=11, r=8, L=13, q=4, \rho=0$)
%
%\bigskip

\item 
$\{0,19, 29, 31, 38, 48, 50, 57, 58, 60, 62, 67, 69, 76,\dots\}$
  \begin{itemize}
  \item  $g=63$,
  \item $c=76$,
  \item $m=19$,
  \item $p=12$,
  \item $r=9$,
  \item $L=13$,
  \item $q=4$,
  \item $\rho=0$.
\end{itemize}

\item
$\{0,19, 30, 31, 38, 49, 50, 57, 60, 61, 62, 68, 69, 76,\dots\}$
\begin{itemize}
  \item $g=63$, \item $c=76$, \item $m=19$, \item $p=12$, \item $r=9$, \item $L=13$, \item $q=4$, \item $\rho=0$.
  \end{itemize}
\end{itemize}

%Our algorithm only checks the Eliahou inequality for semigroups with $l\geq 3$. But as proved in \cite[Section 7.2]{eliahou} this covers all possible counterexamples.

This allows us to state the next result.

\begin{lemma}
The Wilf conjecture holds for all semigroups of genus up to $65$.
\end{lemma}

We notice that for all Eliahou counterexamples of genus up to $65$, $l=13$, $q=4$, $\rho=0$, $g=3$ modulo $4$.

Manuel Delgado constructed in \cite{Delgado}, for each integer number, infinite families of numerical semigroups having Eliahou constant equal to that number.
In particular, he constructed infinite families of semigroups with negative Eliahou constant.
The semigroups in these families are of the form
$$S^{(i, j)} ( p, \tau ) = \langle m^{(i, j)} , g^{(i, j)} , g^{(i, j)} + 1\rangle \mid_{c^{(i, j)}},$$
for $p$ an even positive integer and $\tau,i,j$ non-negative integers, where
\begin{eqnarray*}
  m^{(i, j)} %&=& \frac{p^2}{4}+2p+2 +\tau\frac{p}{2}   +j\frac{p}{2}\\
  &=& \frac{p^2}{4}+p(\frac{\tau}{2}+2)+2+j\frac{p}{2} \\
  g^{(i, j)} &=& 
  %  \frac{p^2}{2}+4p+4+\tau p  -(\tau+1)-\left(\frac{p}{2}+2+1\right)  + j(p-1)+ i m^{(i, j)}\\
  %  \frac{p^2}{2}+4p+4+\tau p  -\frac{p}{2}-4-\tau  + j(p-1)+ i m^{(i, j)}\\
  \frac{p^2}{2}+p(\tau+\frac{7}{2})-\tau  + j(p-1)+ i m^{(i, j)}\\
    c^{(i, j)} &=&
  \frac{p^3}{4}+p^2(\frac{\tau}{2}+2)+2p -\tau + j\frac{p^2}{2}+i\left(\frac{p}{2}+1\right)m^{(i,j)}
\end{eqnarray*}

%We will call these semigroups {\em Delgado semigroups}.
One can check that none of the two semigroups listed above is of this kind.
In the first case it is easy to see, since the difference between the second and third generator is not one.
In the second case, we use that $\rho=\tau$ \cite[Lemma 5.4.8]{Delgado}
and see, by exhaustive search that there is no combination of a positive even integer $p$ and non-negative integers $i,j$ such that
\begin{eqnarray*}
  19 &=& \frac{p^2}{4}+2p+2+j\frac{p}{2} \\
  30 &=& 
  %  \frac{p^2}{2}+4p+4+\tau p  -(\tau+1)-\left(\frac{p}{2}+2+1\right)  + j(p-1)+ i m^{(i, j)}\\
  %  \frac{p^2}{2}+4p+4+\tau p  -\frac{p}{2}-4-\tau  + j(p-1)+ i m^{(i, j)}\\
  \frac{p^2}{2}+\frac{7}{2}p  + j(p-1)+ 19i \\
  76 &=&
  \frac{p^3}{4}+p^2(\frac{\tau}{2}+2)+2p + j\frac{p^2}{2}+19i\left(\frac{p}{2}+1\right)
\end{eqnarray*}

\section*{Acknowledgment}

The authors would like to thank Enric Pons Montserrat, Manuel Delgado, and Julio Fernández-González for their contribution in this work. 
This work was partly supported by the Catalan Government under grant 2017 SGR 00705 and by the Spanish Ministry of Economy and Competitivity under grant \mbox{TIN2016-80250-R} and grant \mbox{RTI2018-095094-B-C21}.

%\bibliographystyle{plain}
%\bibliography{bib}

\begin{thebibliography}{1}

\bibitem{seeds}
M.~Bras-Amor{\'o}s and J.~Fern{\'a}ndez-Gonz{\'a}lez.
\newblock Computation of numerical semigroups by means of seeds.
\newblock {\em Math. Comp.}, 87(313):2539--2550, 2018.

\bibitem{Delgado}
Manuel Delgado.
\newblock On a question of {E}liahou and a conjecture of {W}ilf.
\newblock {\em Math. Z.}, 288(1-2):595--627, 2018.

\bibitem{eliahou}
Shalom Eliahou.
\newblock Wilf's conjecture and {M}acaulay's theorem.
\newblock {\em J. Eur. Math. Soc. (JEMS)}, 20(9):2105--2129, 2018.

\bibitem{FromentinHivert}
J.~Fromentin and F.~Hivert.
\newblock Exploring the tree of numerical semigroups.
\newblock {\em Math. Comp.}, 85(301):2553--2568, 2016.

\bibitem{RG}
J.~C. Rosales and P.~A. Garc\'{\i}a-S\'{a}nchez.
\newblock {\em Numerical semigroups}, volume~20 of {\em Developments in
  Mathematics}.
\newblock Springer, New York, 2009.

\bibitem{wilf}
Herbert~S. Wilf.
\newblock A circle-of-lights algorithm for the ``money-changing problem''.
\newblock {\em Amer. Math. Monthly}, 85(7):562--565, 1978.

\end{thebibliography}

\end{document}